\documentclass[12pt,leqno]{article}
\usepackage{amsmath,amsfonts}
\usepackage[dvips]{graphics,color}
\usepackage{latexsym}
\usepackage{amssymb}
\usepackage{graphicx}
\def\e0{\epsilon}

\def\f0{\phi}

\def\f{\varphi}

\def\e{\varepsilon}

\def\bfabs{\par\bigskip \noindent}

\date{ }
\begin{document}
\title{ Embedded hypersurfaces with constant $m^{\text{th}}$ mean curvature in a
unit  sphere
}
\author{Qing-Ming Cheng, \ Haizhong Li\  and\ Guoxin Wei} \maketitle

\begin{abstract}
\noindent In this paper, we study $n$-dimensional hypersurfaces with
constant $m^{\text{th}}$ mean curvature in a unit sphere $S^{n+1}(1)$
and construct many
compact nontrivial  embedded  hypersurfaces with
constant $m^{\text{th}}$ mean curvature $H_m>0$ in  $S^{n+1}(1)$, for $1\leq m\leq n-1$.
In particular, if the $4^{\text{th}}$ mean curvature $H_4$ takes value between
$\dfrac{1}{(\tan \frac{\pi}{k})^4}$ and $\dfrac{k^4-4}{n(n-4)}$ for
any integer $k\geq3$, then there exists an $n$-dimensional ($n\geq
5$) compact nontrivial embedded hypersurface with constant $H_4$ in
$S^{n+1}(1)$.

\end{abstract}



\section *{1. Introduction}

It is well known that  Alexandrov [1]  and  Montiel-Ros
$\cite{[MR]}$ proved that the standard round spheres are the only
possible oriented compact embedded hypersurfaces with constant
$m^{\text{th}}$ mean curvature $H_m$ in a Euclidean space
$\mathbb{R}^{n+1}$ , for $m\geq 1$. On the other hand,  one knows
that standard round spheres and Clifford hypersurfaces
$S^{l}(a)\times S^{n-l}(b)$, $1\leq l\leq n-1$ are compact embedded
hypersurfaces in a unit sphere $S^{n+1}(1)$. Hence, it is natural to
ask the following:

\bigskip \noindent
{\bf Question:} Do  there exist  compact  embedded hypersurfaces with constant
$m^{\text{th}}$ mean curvature $H_m$ in  $S^{n+1}(1)$ other than the
standard round spheres and Clifford
hypersurfaces?\\

When $m=1$, namely,  when the mean curvature is constant,
Brito-Leite \cite{[BL]} and Perdomo \cite{[P]} have proved  that there exist
 compact embedded hypersurfaces with constant mean curvature $H$ in $S^{n+1}(1)$,
 which are not isometric to the standard round spheres and the Clifford hypersurfaces.

For $m=2$, that is, when the scalar curvature is constant,
Leite  \cite{[Le]} has proved that there exist
compact nontrivial embedded hypersurfaces with constant scalar
curvature $R$ satisfying $(n-1)(n-2)<R<n(n-1)$ in $S^{n+1}(1)$.
Furthermore, Li-Wei \cite{[LW1]} have proved that there exist many
compact nontrivial
embedded hypersurfaces with constant scalar curvature $R$ satisfying
$R>n(n-1)$ in $S^{n+1}(1)$, recently.
But for  $m>2$, one  knows little about
existence of  compact embedded hypersurfaces with constant $m^{\text{th}}$
mean curvature $H_m$ in $S^{n+1}(1)$. In this paper, we  prove that there exist many compact
nontrivial embedded hypersurfaces with constant  $m^{\text{th}}$
mean curvature $H_m>0$  in $S^{n+1}(1)$, for  $1\leq m\leq n-1$. In particular, for $m=4$, we prove
that there exist a lot of compact embedded hypersurfaces with
constant $4^{\text{th}}$ mean curvature $H_4$ in $S^{n+1}(1)$
if it takes value between
$\dfrac{1}{(\tan \frac{\pi}{k})^4}$ and $\dfrac{k^4-4}{n(n-4)}$ for
any integer $k\geq3$.
Furthermore, for $m=1$, our results reduce to the conclusion of Brito-Leite \cite{[BL]} and Perdomo \cite{[P]}.
For $m=2$, we prove that there are many new  compact embedded hypersurfaces with constant
scalar curvature satisfying  $R>n(n-1)$ in $S^{n+1}(1)$, other than
ones of   Li-Wei \cite{[LW1]}.

\section*{2. Preliminaries}

Let $M$ be an $n$-dimensional hypersurface of a unit sphere
$S^{n+1}(1)$ with constant $m^{\text{th}}$ mean curvature $H_m$. We
choose a local orthonormal frame $\{e_{A}\}_{1\leq A\leq n+1}$ in
$S^{n+1}$, with dual coframe $\{\omega_{A}\}_{1\leq A\leq n+1}$,
such that, at each point of $M$, $e_1, \cdots, e_n$ are tangent to
$M$ and $e_{n+1}$ is the positively oriented unit normal vector. We
shall make use of the following convention on the ranges of indices:
$$1\leq A, B, C,\cdots,\leq n+1;\ \ \ \ 1\leq i, j, k, \cdots,\leq
n.$$

Then the structure equations of $S^{n+1}$ are given by
$$d\omega_{A}=\sum_{B=1}^{n+1}\omega_{AB}\wedge \omega_{B}, \ \ \
\omega_{AB}+\omega_{BA}=0,\eqno{(2.1)}$$
$$d\omega_{AB}=\sum_{C=1}^{n+1}\omega_{AC}\wedge\omega_{CB}-\omega_A\wedge\omega_B.\eqno{(2.2)}$$

When restricted to $M$, we have $\omega_{n+1}=0$ and
$$0=d\omega_{n+1}=\sum_{i=1}^{n}\omega_{n+1i}\wedge\omega_i.\eqno{(2.3)}$$

By Cartan's lemma, there exist functions $h_{ij}$ such that
$$\omega_{i n+1}=\sum_{j=1}^{n}h_{ij}\omega_{j}, \ \ \
h_{ij}=h_{ji}.\eqno{(2.4)}$$

This gives the second fundamental form of $M$,
$B=\sum_{i,j}h_{ij}\omega_i\omega_je_{n+1}$. The mean curvature $H$
is defined by $H=\frac{1}{n}\sum_{i}h_{ii}$. From (2.1)-(2.4) we
obtain the structure equations of $M$ (see [8])
$$d\omega_{i}=\sum_{j=1}^{n}\omega_{ij}\wedge \omega_{j}, \ \ \
\omega_{ij}+\omega_{ji}=0,\eqno{(2.5)}$$
$$d\omega_{ij}=\sum_{k=1}^{n}\omega_{ik}\wedge\omega_{kj}-\frac{1}{2}\sum_{k,l=1}^{n}R_{ijkl}\omega_k\wedge\omega_l.
\eqno{(2.6)}$$ and the Gauss equations
$$R_{ijkl}=\delta_{ik}\delta_{jl}-\delta_{il}\delta_{jk}+(h_{ik}h_{jl}-h_{il}h_{jk}).\eqno{(2.7)}$$
$$R-n(n-1)=n(n-1)(r-1)=n^2H^2-S.\eqno{(2.8)}$$
where $R_{ijkl}$ denotes the components of the Riemannian curvature
tensor of $M$, $R=n(n-1)r$ is the scalar curvature of $M$ and
$S=\sum_{i,j=1}^{n}h_{ij}^2$ is the square norm of the second
fundamental form of $M$.

Let $h_{ijk}$ denote the covariant derivative of $h_{ij}$. We then
have (see [8])
$$\sum_{k}h_{ijk}\omega_k=dh_{ij}+\sum_kh_{kj}\omega_{ki}+\sum_kh_{ik}\omega_{kj}.\eqno{(2.9)}$$

Thus, by exterior differentiation of (2.4), we obtain the Codazzi
equation (see [8])
$$h_{ijk}=h_{ikj}.\eqno{(2.10)}$$

We choose $e_1,\cdots,e_n$ such that
$$h_{ij}=\lambda_i\delta_{ij}.\eqno{(2.11)}$$

Let $H_m$ be $m^{\text{th}}$ mean curvature of $M$, then we have
$$C_{n}^{m}H_m=\sum_{1\leq i_1<i_2<\cdots<i_m\leq
n}\lambda_{i_1}\cdots\lambda_{i_m},\eqno{(2.12)}$$ where
$C_{n}^m=\frac{n!}{m!(n-m)!}$.

In \cite{[O]}, Otsuki proved the following

\par\bigskip \noindent{\bf Lemma 2.1} (${\cite{[O]}}$).
{\it Let $M$ be an $n$-dimensional hypersurface in a unit sphere
$S^{n+1}$ such that the multiplicities of principal curvatures are
all constant. Then the distribution of the space of principal
vectors corresponding to each principal curvature is completely
integrable. In particular, if the multiplicity of a principal
curvature is greater than 1, then this principal curvature is
constant on each integral submanifold of the corresponding
distribution of the space of principal vectors. }\\

From Lemma 2.1, we can easily obtain the following theorem.

\bfabs {\bf Theorem 2.1}. {\it Let $M$ be an $n$-dimensional
oriented complete hypersurface in a unit sphere $S^{n+1}$ with
constant $m^{\text{th}}$ mean curvature $H_m$ and with two distinct
principal curvatures. If the multiplicities of these two distinct
principal curvatures are greater than 1, then $M$ is isometric to
Riemannian product $S^{k}(a)\times S^{n-k}(b)$, $2\leq k\leq n-2$.}

\section*{3. A representation formula of principal curvatures}

Now, let us consider that $M$ is an $n$-dimensional oriented
hypersurface with constant $m^{\text{th}}$ mean curvature $H_m$ and
with two distinct principal curvatures in $S^{n+1}$. If
multiplicities of these two distinct principal curvatures are all
great than 1, then we can deduce from Theorem 2.1 that $M$ is
isometric to $S^{k}(a)\times S^{n-k}(b)$, $2\leq k\leq n-2$. Hence,
we shall assume that one of these two distinct principal curvatures
is simple, that is, we assume
$$\lambda_1=\lambda_2=\cdots=\lambda_{n-1}=\lambda,\ \ \
\lambda_n=\mu.\eqno{(3.1)}$$

Since $H_m$ is constant, we obtain from (2.12) that
$$
C_n^mH_m=C_{n-1}^m\lambda^m+C_{n-1}^{m-1}\lambda^{m-1}\mu.\eqno{(3.2)}$$

By Lemma 2.1, let us denote the integral submanifold through $x\in
M$, corresponding to $\lambda$ by $M_1^{n-1}(x)$. We write
$$d\lambda=\sum_{i}\lambda_{,i}\omega_i,\ \ \
d\mu=\sum_{j}\mu_{,j}\omega_j.\eqno{(3.3)}$$

We assume that $\lambda>0$ on $M$. Then Lemma 2.1 implies
$$\lambda_{,1}=\cdots=\lambda_{,n-1}=0.\eqno{(3.4)}$$

Then (3.2) yields
$$\mu=\frac{C_n^mH_m-C_{n-1}^m\lambda^m}{C_{n-1}^{m-1}\lambda^{m-1}}=\frac{nH_m-(n-m)\lambda^m}{m\lambda^{m-1}},
\eqno{(3.5)}$$ and from the formula
$$\lambda-\mu=\frac{n(\lambda^m-H_m)}{m\lambda^{m-1}},\eqno{(3.6)}$$
we obtain that $$\lambda^m-H_m\neq0.\eqno{(3.7)}$$

By means of (2.9) and (2.11), we obtain
$$\sum_kh_{ijk}\omega_k=\delta_{ij}d\lambda_j+(\lambda_i-\lambda_j)\omega_{ij}.\eqno{(3.8)}$$

We adopt the notational convention that
$$1\leq a, b, c, \cdots\leq n-1.$$

 From (3.1), (3.2) and (3.8), we have
$$h_{ijk}=0,\ \ \  \mbox{if}\ i\neq j,\ \lambda_i=\lambda_j,\eqno{(3.9)}$$
$$h_{aab}=0,\ \ \ h_{aan}=\lambda_{,n},\eqno{(3.10)}$$
$$h_{nna}=0,\ \ \ h_{nnn}=\mu_{,n}. \eqno{(3.11)}$$

Combining this with (2.10) and the formula
$$\sum_{i}h_{ani}\omega_i=dh_{an}+\sum_ih_{in}\omega_{ia}+\sum_ih_{ai}\omega_{in}=(\lambda-\mu)\omega_{an},\eqno{(3.12)}
$$
we obtain from (3.10) and (3.6)
$$\omega_{an}=\frac{\lambda_{,n}}{\lambda-\mu}\omega_a=\frac{m\lambda^{m-1}\lambda_{,n}}{n(\lambda^m-H_m)}\omega_a.
\eqno{(3.13)}$$

Therefore we have
$$d\omega_n=\sum_{a}\omega_{na}\wedge\omega_a=0.\eqno{(3.14)}$$

Notice that we may consider $\lambda$ to be locally a function of
the parameter $s$, where $s$ is the arc length of an orthogonal
trajectory of the family of the integral submanifolds corresponding
to $\lambda$. We may put
$$\omega_n=ds.$$

Thus, for $\lambda=\lambda(s),$ we have
$$d\lambda=\lambda_{,n}ds,\ \ \ \lambda_{,n}=\lambda^{'}(s).\eqno{(3.15)}$$

From (3.6) and (3.13), we get
$$\omega_{an}=\frac{m\lambda^{m-1}\lambda_{,n}}{n(\lambda^m-H_m)}\omega_a=
\frac{m\lambda^{m-1}\lambda^{'}(s)}{n(\lambda^m-H_m)}\omega_a=\{\log\mid
\lambda^m-H_m\mid ^{1/n}\}^{'}\omega_a,\eqno{(3.16)}$$ which shows
that the integral submanifolds $M_1^{n-1}(x)$ corresponding to
$\lambda$ is umbilical in $M$ and $S^{n+1}$.

On the other hand, we can deduce from (3.16) that
$$\nabla_{e_n}e_{n}=\sum_{k=1}^{n}\omega_{ni}(e_n)e_i=0.$$

According to the definition of geodesic, we know that the integral
curve of the principal vector field $e_n$ corresponding to the
principal curvature $\mu$ is a geodesic.\\

This proves the following result:
\par\bigskip \noindent{\bf Lemma 3.1}. {\it If $M$ is an $n$-dimensional oriented complete hypersurface $(n\geq3)$ in
$S^{n+1}$ with constant $m^{\text{th}}$ mean curvature $H_m$ and
with two distinct principal curvatures,
one of which is simple, then \\
$(1)$ the integral submanifold $M_{1}^{n-1}(x)$ through $x\in M$
corresponding to $\lambda$ is umbilical in $M$ and
$S^{n+1}$,\\
$(2)$ the integral curve of the principal vector field $e_n$
corresponding to the principal curvature $\mu$ is a geodesic.}\\

 Now we state our Theorem 3.1 as follows:

\par\bigskip \noindent{\bf Theorem 3.1}. {\it If $M$ is an $n$-dimensional oriented complete hypersurface $(n\geq3)$ in
$S^{n+1}$ with constant $m^{\text{th}}$ mean curvature $H_m$ and
with two distinct principal curvatures one of which is simple, then
$M$ is isometric to a complete hypersurface of $S^{n-1}(c(s))\times
M^1$, where $S^{n-1}(c(s))$ is of constant curvature
$[(\log\mid\lambda^m-H_m\mid^{1/n})^{'}]^2+\lambda^2+1$. And
$w=\mid\lambda^m-H_m\mid^{-1/n}$ satisfies the following ordinary
differential equation of order 2:
$$\frac{d^2w}{ds^2}-w\left\{\frac{(n-m)(w^{-n}+H_m)^{(2-m)/m}}{mw^n}-H_m(w^{-n}+H_m)^
{(2-m)/m}-1 \right\}=0 .\eqno{(3.17)}$$}

\vskip 3pt\noindent {\it The proof of  Theorem 3.1}. According to
the structure equations of $S^{n+1}$ and (3.16), we may compute
\begin{align}
d\omega_{an}&=\sum_{b=1}^{n-1}\omega_{ab}\wedge\omega_{bn}+\omega_{an+1}\wedge\omega_{n+1n}-\omega_a\wedge \omega_n\notag\\
            &=(\log\mid\lambda^m-H_m\mid
            ^{1/n})^{'}\sum_{b=1}^{n-1}\omega_{ab}\wedge\omega_b-\lambda\mu\omega_a\wedge
            ds-\omega_a\wedge ds,\notag
 \end{align}
\begin{align}
d\omega_{an}&=d[(\log\mid \lambda^m-H_m\mid^{1/n})^{'}\omega_a]\notag\\
            &=\{\log\mid\lambda^m-H_m\mid^{1/n}\}^{''}ds\wedge\omega_a+\{\log\mid
            \lambda^m-H_m\mid^{1/n}\}^{'}d\omega_a\notag\\
            &=\left\{-(\log\mid\lambda^m-H_m\mid^{1/n})^{''}+\left[(\log\mid\lambda^m-H_m\mid^{1/n})^{'}\right]^2\right
            \}\omega_a\wedge
            ds\notag\\
            &\ \ \
            +\left(\log\mid\lambda^m-H_m\mid^{1/n}\right)^{'}\sum_{b=1}^{n-1}\omega_{ab}\wedge\omega_{b}.\notag
 \end{align}

Then we obtain from two equalities above that
$$\left\{\log\mid\lambda^m-H_m\mid^{1/n}\right\}^{''}-\left[(\log\mid\lambda^m-H_m\mid^{1/n})^{'}\right]^2
-\lambda\mu-1=0.\eqno{(3.18)}$$

Combining (3.18) with (3.6), we have
$$\left\{\log\mid\lambda^m-H_m\mid^{1/n}\right\}^{''}-\left[(\log\mid\lambda^m-H_m\mid^{1/n})^{'}\right]^2
+\frac{(n-m)\lambda^m-nH_m }{m\lambda^{m-2}}-1=0.\eqno{(3.19)}$$

We know that $\lambda^m-H_m\neq0$. If $\lambda^m-H_m<0$, from (3.6),
we have
$$\lambda^2-\lambda\mu=\frac{n(\lambda^m-H_m)}{m\lambda^{m-2}}<0,\
\ \ H_m>0.\eqno{(3.20)}$$
According to the Gauss equation (2.7), we
know that the sectional curvature of $M$ is not less than $1$ and
$H_m>0$. By a direct calculation, we know that $M$ is isometric to a
totally umbilical hypersurface. This is impossible because $M$ has
two distinct principal curvatures. Hence, $\lambda^m-H_m>0$. Let us
define a positive function $w(s)$ over $s\in(-\infty, +\infty)$ by
$$w=(\lambda^m-H_m)^{-1/n},\eqno{(3.21)}$$
then (3.19) reduces to
$$\frac{d^2w}{ds^2}-w\left\{\frac{(n-m)(w^{-n}+H_m)^{(2-m)/m}}{mw^n}-H_m(w^{-n}+H_m)^
{(2-m)/m}-1 \right\}=0.\eqno{(3.22)}$$

Integrating (3.22), we obtain
$$\left(\frac{dw}{ds}\right)^2=C-w^2\left(w^{-n}+H_m\right)^{\frac{2}{m}}-w^2,\eqno{(3.23)}$$
where $C$ is the constant of integration.

We consider the frame $\{x,e_1,e_2,\cdots,e_n,e_{n+1}\}$ in the
Euclidean space $\mathbb{R}^{n+2}$. Then, by (2.4), (3.13) and
(3.18), we obtain
\begin{align}
de_a&=\sum_{b=1}^{n-1}\omega_{ab}e_b+\omega_{an}e_n+\omega_{an+1}e_{n+1}-\omega_a e_{n+2}\notag\\
    &=\sum_{b=1}^{n-1}\omega_{ab}e_b+(\log\mid\lambda^m-H_m\mid ^{1/n})^{'}\omega_ae_n-\lambda\omega_a
    e_{n+1}-\omega_a e_{n+2}\notag\\
    &=\sum_{b=1}^{n-1}\omega_{ab}e_b+\left\{(\log\mid\lambda^m-H_m\mid ^{1/n})^{'}-\lambda
    e_{n+1}-e_{n+2}\right\}\omega_a\notag
 \end{align}
\begin{align}
d&\left\{(\log\mid\lambda^m-H_m\mid ^{1/n})^{'}-\lambda
    e_{n+1}-e_{n+2}\right\}\ \ \ \ \ \ \ \ \ \ \ \ \ \ \ \ \notag\\
    &\ \ \ \ \ \ =\left\{(\log\mid\lambda^m-H_m\mid ^{1/n})^{''}-\lambda^{'}
    e_{n+1}\right\}ds\notag\\
    &\ \ \ \ \ \ \ \ +(\log\mid\lambda^m-H_m\mid ^{1/n})^{'}(\sum_{a=1}^{n-1}\omega_{na}e_a
    +\omega_{nn+1}e_{n+1})-\lambda\left(\sum_{a=1}^{n-1}\omega_{n+1a}e_a+\omega_{n+1n}e_n\right)\notag\\
    &\ \ \ \ \ \ \ \ -\sum_{a=1}^{n-1}\omega_a e_a-\omega_n e_n\notag\\
    &\ \ \ \ \ \ \equiv \left\{(\log\mid\lambda^m-H_m\mid
    ^{1/n})^{''}-\lambda\mu-1\right\}e_n\omega_n\notag\\
   &\ \ \ \ \ \ \ \ -\left\{\lambda^{'}+(\log\mid\lambda^m-H_m\mid
   ^{1/n})^{'}\mu\right\}e_{n+1}\omega_n\ \ (\mbox{mod}\
   \{e_1,\cdots,e_{n-1}\})\notag\\
    &\ \ \ \ \ \ =(\log\mid\lambda^m-H_m\mid ^{1/n})^{'}\left\{(\log\mid\lambda^m-H_m\mid
  ^{1/n})^{'}e_n-\lambda e_{n+1}-e_{n+2}\right\}ds\notag
\end{align}

By putting
$$W=e_1\wedge e_2\wedge\cdots\wedge e_{n-1}\wedge\left\{(\log\mid\lambda^m-H_m\mid
  ^{1/n})^{'}e_n-\lambda e_{n+1}-e_{n+2}\right\},\eqno{(3.24)}$$
we can show that
$$dW=(\log\mid\lambda^m-H_m\mid ^{1/n})^{'}Wds.\eqno{(3.25)}$$

(3.25) shows that $n$-vector $W$ in $\mathbb{R}^{n+2}$ is constant
along $M_1^{n-1}(x)$. Hence there exists an $n$-dimensional linear
subspace $E^n(s)$ in $\mathbb{R}^{n+2}$ containing $M_1^{n-1}(x)$.
(3.25) also implies that the $n$-vector field $W$ only depends on
$s$ and by integrating it, we get
$$W=\left\{\frac{\lambda^m(s)-H_m}{\lambda^m(s_0)-H_m}\right\}^{1/n}W(s_0).\eqno{(3.26)}$$
Theorefore, we have that $E^n(s)$ is parallel to $E^n(s_0)$ in
$\mathbb{R}^{n+2}$ for every $s$.

From the calculation
\begin{align}
d\omega_{ab}-\sum_{c=1}^{n-1}\omega_{ac}\wedge\omega_{cb}&=\omega_{an}\wedge\omega_{nb}+\omega_{an+1}\wedge
                                                                  \omega_{n+1b}-\omega_a\wedge\omega_b\notag\\
                       &=-\left\{[(\log\mid\lambda^m-H_m\mid
                       ^{1/n})^{'}]^2+\lambda^2+1\right\}\omega_a\wedge\omega_b,\notag
 \end{align}
we see that the curvature of $M_1^{n-1}(x)$ is
$[(\log\mid\lambda^m-H_m\mid^{1/n})^{'}]^2+\lambda^2+1$ and
$M_1^{n-1}(x)$ is locally isometric to $S^{n-1}(c(s))$. Therefore,
$M$ is isometric to a complete hypersurface of revolution
$S^{n-1}(c(s))\times M^1$.

 This proves Theorem 3.1.

\section*{4.  A representation formula of radius }

One knows that the following immersion:
$$x: M^n\hookrightarrow S^{n+1}(1)\subset R^{n+2},$$
$$(s, t_1, \cdots,  t_{n-1})\mapsto
(y_1(s)\varphi_1,\cdots,y_1(s)\varphi_n,y_{n+1}(s),y_{n+2}(s)).\eqno{(4.1)}$$
$$\varphi_i=\varphi_i(t_1,\cdots,t_{n-1}),\ \ \
\varphi_1^2+\cdots+\varphi_n^2=1\eqno{(4.2)}$$ is a parametrization
of a rotational hypersurface generated by a curve $(y_1(s),
y_{n+1}(s),y_{n+2}(s))$. Since the curve $(y_1(s), y_{n+1}(s),
y_{n+2}(s))$ belongs to $S^{2}(1)$ and the parameter $s$ can be
chosen as its arc length, we have
$$y_1^2(s)+y_{n+1}^2(s)+y_{n+2}^2(s)=1,\ \ \
\dot{y}_1^2(s)+\dot{y}_{n+1}^2(s)+\dot{y}_{n+2}^2(s)=1\eqno{(4.3)}$$
where the dot denotes the derivative with respect to $s$ and from
(4.3) we can obtain $y_{n+1}(s)$ and $y_{n+2}(s)$ as functions of
$y_1(s)$. In fact, we can write
$$y_1(s)=\cos \vartheta(s),\ \ y_{n+1}(s)=\sin \vartheta(s)\cos \theta(s),\ \
y_{n+2}(s)=\sin \vartheta(s)\sin \theta(s).\eqno{(4.4)}$$ We can
deduce from (4.3) that
$$\dot{\vartheta}^2+\dot{\theta}^2\sin^2\vartheta=1.\eqno{(4.5)}$$

It follows from equation (4.5) that $\dot{\vartheta}^2\leq1$.
Combining these with
$\dot{\vartheta}^2=\frac{\dot{y}_1^2}{1-y_1^2}$, we have
$$\dot{y}_1^2+y_1^2\leq1.\eqno{(4.6)}$$

We can get the plane curve $\zeta$ from $\alpha$ by projection of
$S_{+}^2=\{(y_1,y_{n+1},y_{n+2})\mid y_1\geq 0,
y_1^2+y_{n+1}^2+y_{n+2}^2=1\}$ onto the unit disk
$E=\{(y_{n+1},y_{n+2})\mid y_{n+1}^2+y_{n+2}^2\leq1\}$. Then the
plane curve $\zeta$ can be written as
$$y_{n+1}(s)=\sin \vartheta(s)\cos \theta(s),\ \
y_{n+2}(s)=\sin \vartheta(s)\sin \theta(s).\eqno{(4.7)}$$

Writing $r(s)=y_1(s)$, (4.5) can be written as

$$\dot{\theta}^2=\frac{1-\dot{\vartheta}^2}{\sin^2\vartheta}=\frac{1-r^2-\dot{r}^2}{(1-r^2)^2}.\eqno{(4.8)}$$

Do Carmo and Dajczer proved the following

\par\bigskip \noindent{\bf Lemma 4.1} (\cite{[DD]}). {\it Let
$M^n$ be a rotational hypersurface of $S^{n+1}(1)$. Then the
principal curvatures $\lambda_i$ of $M^n$ are
$$\lambda_i=\lambda=-\frac{\sqrt{1-r^2-\dot{r}^2}}{r}\eqno{(4.9)}$$
for $i=1,\cdots,n-1$, and
$$\lambda_n=\mu=\frac{\ddot{r}+r}{\sqrt{1-r^2-\dot{r}^2}}.\eqno{(4.10)}$$}

On the other hand, let us fix a point $p_0\in M$, let $\gamma(u)$ be
the only geodesic in $M$ such that $\gamma(0)=p_0$ and
$\gamma^{'}(0)=e_n(p_0)$. From (3.16) of Section 3, we know that
$\gamma(u)=e_n(\gamma(u))$. Note that $\gamma(u)$ is also a line of
curvature. Let us denote by $g(u)=w(\gamma(u))$. Since $H_m$ is
constant, we know from (3.23) that

$$(g^{'})^2+g^2\left(g^{-n}+H_m\right)^{\frac{2}{m}}+g^2=C.\eqno{(4.11)}$$

From (4.11), we have $C>0$. Moreover, by a direct calculation, we
get

$$q(x)=C-x^2\left(x^{-n}+H_m\right)^{\frac{2}{m}}-x^2\eqno{(4.12)}$$
is positive on a interval $(t_1, t_2)$ with $0<t_1<t_2$ and
$q(t_1)=q(t_2)=0$. From (4.11), we know that $g(u)$ is periodic. And
the period is the following

$$T=2\int_{t_1}^{t_2}\frac{dt}{\sqrt{C-t^2\left(t^{-n}+H_m\right)^{\frac{2}{m}}-t^2}}dt.\eqno{(4.13)}$$

From (4.1) and Theorem 3.1, we have
$$\frac{1}{r^2}=[(\log\mid\lambda^m-H_m\mid^{1/n})^{'}]^2+\lambda^2+1.$$

Then we know form (3.23), (4.11) that

$$r(u)=\frac{g(u)}{\sqrt{C}},\ \ \ \ \
g(u)=(\lambda^m-H_m)^{-\frac{1}{n}}.\eqno{(4.14)}$$

From (4.8), (4.9) and (4.11), we obtain the period $P(H_m,n,c)$ of
hypersurfaces

\begin{equation*}
\aligned &\ \ \ \ P(H_m,n,C)\\
&=\theta(T)=\int_0^{T}\frac{\sqrt{1-r^2-\dot{r}^2}}{1-r^2}ds\\
&=\int_0^{T}\frac{r(s)\lambda(s)}{1-r^2(s)}ds\\
&=2\int_0^{\frac{T}{2}}\frac{r(s)\lambda(s)}{1-r^2(s)}ds.
\endaligned
\eqno{(4.15)}
\end{equation*}

\section*{5. Embedded hypersurfaces with constant $H_m>0$}

At first, we give the following Lemma due to Perdomo \cite{[P]}

\par\bigskip \noindent{\bf Lemma 5.1}.{\it Let $\epsilon$ and $\delta$ be positive numbers and
$f: (t_0-\epsilon, t_0+\epsilon)\rightarrow R$ and $y:(-\delta,
\delta)\times(t_0-\epsilon, t_0+\epsilon)\rightarrow R$ be two
smooth functions such that $f(t_0)=f^{'}(t_0)=0$ and
$f^{''}(t_0)=-2a<0$. If for any small $c>0$, $t_1(c)<t_0<t_2(c)$ are
such that $f(t_1(c))+c=0=f(t_2(c))+c$, then $$\lim_{c\rightarrow
0^{+}}\int_{t_1(c)}^{t_{2(c)}}\frac{y(c,t)dt}{\sqrt{f(t)+c}}=\frac{y(0,t_0)\pi}{\sqrt{a}}.$$}

Now we state our  main theorem.

\par\bigskip \noindent{\bf Theorem 5.1}. {\it For any $n\geq 5$ and any
integer $k\geq 3$, if $4^{\text{th}}$ mean curvature $H_4$ takes
value between $\dfrac{1}{(\tan \frac{\pi}{k})^4}$ and
$\dfrac{k^4-4}{n(n-4)}$, then there exists an $n$-dimensional
compact nontrivial embedded hypersurface with constant $H_4>0$ in
$S^{n+1}(1)$.}

 \vskip 3pt\noindent {\it Proof}. Let us rewrite (4.11) as
 $$(g^{'})^2=q(g),\ \ \text{where}\ \
 q(v)=C-v^2(v^{-n}+H_m)^{\frac{2}{m}}-v^2.\eqno{(5.1)}$$

 We know that for some value of $C$, the function $q$ has positive
 values between two positive roots of $q$, denoted by $t_1$ and
 $t_2$. A direct calculation shows that
 $$q^{'}(v)=2v\left\{-(v^{-n}+H_m)^{\frac{2}{m}}+\frac{n}{m}v^{-n}(v^{-n}+H_m)^{\frac{2-m}{m}}-1\right\}.\eqno{(5.2)}$$

\begin{equation*}
\aligned &\ \ \ q^{''}(v)\\
&=-\frac{2(v^{-n}+H_m)^{\frac{2-2m}{m}}}{m^2}((2n^2-3nm+m^2)v^{-2n}+m(n^2-3n+2m)H_mv^{-n}+m^2H_m^2)-2\\
&<-2.
\endaligned
\eqno{(5.3)}
\end{equation*}

If $m=4$ and $H_4=1$, we have the only positive root of $q^{'}$ is
$$v_0=(\frac{(n-4)^2}{8n-16})^{\frac{1}{n}}.\eqno{(5.4)}$$
Therefore, for positive values of $v$, the function $q$ increase
from $0$ to $v_0$ and decrease for values greater than $v_0$.  Then
we obtain that $q(v_0)=C-c_0$, where
$$c_0=v_0^2((v_0^{-n}+1)^{\frac{1}{2}}+1)=(\frac{(n-4)^2}{8n-16})^{\frac{2}{n}}\times
(\frac{n}{n-4}+1).\eqno{(5.5)}$$

Therefore, whenever $C>c_0$, we will have the two positive roots of
the function $q(v)$ that we will denote by $t_1(C)$ and $t_2(C)$. By
computing, we have $q^{''}(v_0)=-2a$, where
$$a=\frac{2(n-2)^2}{n}.\eqno{(5.6)}$$

Hence, we get from (4.15) that
$$P(H_4,n,C)=2\int_0^{\frac{T}{2}}\frac{r(s)\lambda(s)}{1-r^2(s)}ds.\eqno{(5.7)}$$
Since $r(s)=\frac{g(s)}{\sqrt{C}}$ and
$\lambda(s)=(g^{-n}+1)^{\frac{1}{4}}$, we have
$$P(H_4,n,C)=2\int_0^{\frac{T}{2}}\frac{\sqrt{C}g(s)(g^{-n}(s)+1)^{\frac{1}{4}}}{C-g^2(s)}ds.\eqno{(5.8)}$$

Since $g(0)=t_1(C)$ and $g(\frac{T}{2})=t_2(C)$, by doing the
substitutions $t=g(s)$, we have
$$P(H_4,n,C)=2\int_{t_1(C)}^{t_2(C)}\frac{\sqrt{C}t(t^{-n}+1)^{\frac{1}{4}}}{C-t^2}\frac{1}{\sqrt{q(t)}}ds.\eqno{(5.9)}$$

Using Lemma 5.1, we have
$$\lim_{C\rightarrow c_0^{+}}P(H_4=1,n,C)=\frac{2\pi}{\sqrt{a}}
\frac{\sqrt{c_0}v_0(v_0^{-n}+1)^{\frac{1}{4}}}{c_0-v_0^2}=\frac{2\pi\sqrt{n-2}}{n-2}.\eqno{(5.10)}$$

On the other hand, we will estimate $P(H_m,n,C)$ when
$C\rightarrow\infty$, we make the substitution $t=r(s)$ and obtain
$$P(H_m,n,C)=2\int_{\frac{t_1(C)}{\sqrt{C}}}^{\frac{t_2(C)}{\sqrt{C}}}\frac{t((\sqrt{C}t)^{-n}+H_m)^{\frac{1}{m}}}
{(1-t^2)\sqrt{1-t^2(1+(H_m+(\sqrt{C}t)^{-n})^{\frac{2}{m}}}}.\eqno{(5.11)}$$

Since
$$\tilde{q}=1-t^2(1+(H_m+(\sqrt{C}t)^{-n})^{\frac{2}{m}}\eqno{(5.12)}$$ have
two positive roots converge to $0$ and
$\frac{1}{\sqrt{1+H_m^{\frac{2}{m}}}}$, we obtain

$$\lim_{C\rightarrow
\infty}P(H_m,n,C)=2\int_0^{\frac{1}{\sqrt{1+H_m^{\frac{2}{m}}}}}
\frac{tH_m^{\frac{1}{m}}}{(1-t^2)\sqrt{1-t^2(1+H_m^{\frac{2}{m}}})}dt=2\arctan\frac{1}{H_m^{\frac{1}{m}}}.\eqno{(5.13)}$$

If $m=4$ and $H_4=1$, we have that

$$\lim_{C\rightarrow \infty}P(H_4=1,n,C)=2\arctan
\frac{1}{(H_4)^{\frac{1}{4}}}=\frac{\pi}{2}.\eqno{(5.14)}$$\\



Next, we consider the case $m=4$ and $0<H_4\neq 1$.

In this case, we have the only positive root of $q^{'}$ is
$$v_0=(\frac{\sqrt{n(n-4)H_4+4}-nH_4+4H_4-2}{4H_4(1-H_4)})^{\frac{1}{n}}.\eqno{(5.15)}$$

A direct calculation shows that $q(v_0)=C-c_0$, where

\begin{equation*}
\aligned
 c_0
&=v_0^2(v_0^{-n}+H_4)^{\frac{1}{2}}+v_0^2\\
&=(\frac{\sqrt{n(n-4)H_4+4}-nH_4+4H_4-2}{4H_4(1-H_4)})^{\frac{2}{n}}\\
&\ \ \ \times
[(\frac{H_4(\sqrt{n(n-4)H_4+4}-nH_4+2)}{\sqrt{n(n-4)H_4+4}-nH_4+4H_4-2})^{\frac{1}{2}}+1].
\endaligned
\eqno{(5.16)}
\end{equation*}

\begin{equation*}
\aligned &q^{''}(v_0)\\
 &=-2a\\
&=\frac{-2H_4^{\frac{1}{2}}}{(|\sqrt{n(n-4)H_4+4}-nH_4+2|)^{\frac{3}{2}}
(|\sqrt{n(n-4)H_4+4}-nH_4+4H_4-2|)^{\frac{1}{2}}}\\
& \times
\{n^2(n-4)H_4^2+n(-n^2+4n+4)H_4-4n+[n^2-2n+(-n^2+2n)H_4]\sqrt{n(n-4)H_4+4}\}.
\endaligned
\eqno{(5.17)}
\end{equation*}

Therefore, whenever $C>c_0$, we will have two positive roots of the
function $q(v)$ that we will denote by $t_1$ and $t_2$.

Using the results of section 4, we have from (4.15) that
$$P(H_4,n,C)=2\int_0^{\frac{T}{2}}\frac{r(s)\lambda(s)}{1-r^2(s)}ds.\eqno{(5.18)}$$

From (4.14), we have $r(s)=\frac{g(s)}{\sqrt{C}}$ and
$\lambda(s)=(g(s)^{-n}+H_m)^{\frac{1}{m}}$, then we get that
$$P(H_4,n,C)=2\int_0^{\frac{T}{2}}\frac{\sqrt{C}g(s)(g(s)^{-n}+H_4)^{\frac{1}{4}}}{C-g^2(s)}ds.\eqno{(5.19)}$$

Since $g(0)=t_1$ and $g(\frac{T}{2})=t_2$, by doing the substitution
$t=g(s)$ and using Lemma 5.1, we obtain

\begin{equation*}
\aligned &\lim_{C\rightarrow c_0^{+}}
P(H_4,n,C)\\
 &=\frac{2\pi\sqrt{c_0}}{\sqrt{a}\sqrt{c_0-v_0^2}}\\
&=2\pi\frac{|(n-2)(n-nH_4)+(nH_4-n)\sqrt{n(n-4)H_4+4}|^{\frac{1}{2}}}{|n^2(n-4)H_4^2+
n(-n^2+4n+4)H_4-4n+[n^2-2n+(-n^2+2n)H_4]\sqrt{n(n-4)H_4+4}|^{\frac{1}{2}}}\\
&=2\pi\frac{|(n-2)-\sqrt{n(n-4)H_4+4}|^{\frac{1}{2}}}{|(n(4-n)H_4-4)+(n-2)\sqrt{n(n-4)H_4+4}|^{\frac{1}{2}}}\\
&=\frac{2\pi}{[n(n-4)H_4+4]^{\frac{1}{4}}}.
\endaligned
\eqno{(5.20)}
\end{equation*}

On the other hand, we know that
$$\lim_{C\rightarrow
\infty}P(H_4,n,C)=2\arctan\frac{1}{H_4^{\frac{1}{4}}}.\eqno{(5.21)}$$

Therefore, for any fixed $H_4>0$, the function $P(H_4,n,C)$ takes
all the values between
$$A(H_4)=2\arctan\frac{1}{H_4^{\frac{1}{4}}},\ \
B(H_4)=\frac{2\pi}{[n(n-4)H_4+4]^{\frac{1}{4}}}. \eqno{(5.22)}$$

By a direct computation, we know that $A(H_4)$ and $B(H_4)$ are
decreasing functions. Since
$$A(\frac{1}{(\tan
\frac{\pi}{k})^4})=B(\frac{k^4-4}{n(n-4)})=\frac{2\pi}{k},\eqno{(5.23)}$$
where $k\geq3$ is any integer, we deduce that the number
$\frac{2\pi}{k}$ lies between $A(H_4)$ and $B(H_4)$ since they are
decreasing functions, hence, for some constant $C_1$, we have that
$P(H_4,n,C_1)=\frac{2\pi}{k}$. If the period is $\frac{2\pi}{k}$,
then there exists a compact embedded hypersurfaces with constnat
$H_4$ which is not isometric to a round sphere or a Clifford
hypersurface.

We complete the proof of Theorem 5.1.\\

For constant $H_m>0$, we can prove the following

\par\bigskip \noindent{\bf Theorem 5.2}. {\it
For any integer $1\leq m\leq n-1$,  there exist many  nontrivial
embedded hypersurfaces with constant $H_m>0$ in $S^{n+1}(1)$ .}

 \vskip 3pt\noindent {\it Proof}.  By using the similar arguments
 with the proof of Theorem 5.1, we have that

$$v_0=(\frac{n-m}{m})^{\frac{m}{2n}},\ \ c_0=(\frac{n-m}{m})^{\frac{m}{n}}\times \frac{n}{n-m},$$
$$q^{''}(v_0)=-2a=-\frac{4n}{m}.$$

From (4.14), we have $r=\frac{g}{\sqrt{C}}$ and
$\lambda=g^{-\frac{n}{m}}$, then we obtain

$$\lim_{C\rightarrow c_0^{+}}
P(H_m=0,n,C)=\frac{2\pi\sqrt{c_0}}{\sqrt{a}\sqrt{c_0-v_0^2}}
=\sqrt{2}\pi,$$
 by continuity arguments, we can fix
$H_m$ sufficiently small such that
$$\lim_{C\rightarrow c_0^{+}} P(H_m,n,C)>\pi.$$

On the other hand, we deduce from $H_m>0$ that
$$\lim_{C\rightarrow
\infty}P(H_m,n,C)=2\arctan\frac{1}{H_m^{\frac{1}{m}}}<\pi.$$ Hence,
there exists $C_2$, such that $P(H_m,n,C_2)=\pi$.

We complete the proof of Theorem 5.2.

\par\bigskip \noindent{\bf Remark 5.1}. When $m=1$, Theorem 5.2
reduces to the results of Brito and Leite $\cite{[BL]}$.\\

Using the similar arguments as above, we can obtain the following:\\

\noindent When $m=1$, we have

\par\bigskip \noindent{\bf Proposition 5.1} ([14]). {\it For any $n\geq2$ and any
integer $k\geq 2$, if mean curvature $H$ takes value between
$\dfrac{1}{(\tan \frac{\pi}{k})}$ and
$\dfrac{(k^2-2)\sqrt{n-1}}{n\sqrt{k^2-1}}$, then there exists an
$n$-dimensional compact nontrivial embedded hypersurface with
constant mean curvature $H>0$ in $S^{n+1}(1)$.}

\par\bigskip \noindent{\bf Remark 5.2}. Proposition 5.1 is also proved
by Perdomo \cite{[P]}.\\

\noindent When $m=2$, we have

\par\bigskip \noindent{\bf Proposition 5.2}. {\it For any $n\geq3$ and any
integer $k\geq 2$, if $H_2=\frac{R-n(n-1)}{n(n-1)}$ takes value
between $\dfrac{1}{(\tan \frac{\pi}{k})^2}$ and $\dfrac{k^2-2}{n}$,
then there exists an $n$-dimensional compact nontrivial embedded
hypersurface $M$ with constant 2-th mean curvature $H_2>0$ $($i.e.
scalar curvature $R>n(n-1)$$)$ in  $S^{n+1}(1)$, where $R$ is the
scalar curvature of $M$.}
 \vskip 3pt\noindent {\it Proof}.  By using the similar arguments
 with the proof of Theorem 5.1, we have that
 $$v_0=(\frac{n-2}{2(H_2+1)})^{\frac{1}{n}},\ \
 c_0=(\frac{n-2}{2(H_2+1)})^{\frac{2}{n}}\times\frac{n(H_2+1)}{n-2},$$

$$q^{''}(v_0)=-2a=-2n(H_2+1).$$
$$\lim_{C\rightarrow c_0^{+}}
P(H_2,n,C)=\frac{2\pi\sqrt{c_0}}{\sqrt{a}\sqrt{c_0-v_0^2}}=\frac{2\pi}{\sqrt{nH_2+2}}.$$

On the other hand, $$\lim_{C\rightarrow
\infty}P(H_2,n,C)=2\arctan\frac{1}{H_2^{\frac{1}{2}}}.$$

Therefore, for any fixed $H_2>0$, the function $P(H_2,n,C)$ takes
all the values between
$$E(H_2)=2\arctan\frac{1}{H_2^{\frac{1}{2}}},\ \
F(H_2)=\frac{2\pi}{\sqrt{nH_2+2}}. $$

By a direct computation, we know that $E(H_2)$ and $F(H_2)$ are
decreasing functions. Since
$$E(\frac{1}{(\tan
\frac{\pi}{k})^2})=F(\frac{k^2-2}{n})=\frac{2\pi}{k},$$ where
$k\geq2$ is any integer, we deduce that the number $\frac{2\pi}{k}$
lies between $E(H_4)$ and $F(H_4)$ since they are decreasing
functions, hence, for some constant $C_3$, we have that
$P(H_2,n,C_3)=\frac{2\pi}{k}$. If the period is $\frac{2\pi}{k}$,
then there exists a compact embedded hypersurfaces with constnat
$H_2$ $($i.e. constant scalar curvature$)$ which is not isometric to
a round sphere or a Clifford hypersurface.

We complete the proof of Proposition 5.2.

\par\bigskip \noindent{\bf Remark 5.3}. Since
$H_2=\frac{R-n(n-1)}{n(n-1)}$, by a direct calculation, we know that
when $3\leq n\leq 6$, Proposition 5.2 reduces to Theorem 1.1 and
Theorem 1.2 due to Li-Wei \cite{[LW1]}; when $n>6$ and $k=2$,
Proposition 5.2 reduces to Theorem 1.3 due to Li-Wei \cite{[LW1]}.
In Proposition 5.2, we find there exist a lot of new examples
satisfying $R>n(n-1)$. Hence, Proposition 5.2 is the generalization
of  Li-Wei's results \cite{[LW1]}.

\par\bigskip \noindent{\bf Remark 5.4}. For some special $4\neq m>3$, we
 can also obtain some nontrivial embedded hypersurfaces with
$H_m=$constant in $S^{n+1}(1)$ using the same methods.

\begin{flushleft}
\medskip\noindent
\begin{tabbing}
XXXXXXXXXXXXXXXXXXXXXXXXXX*\=\kill Qing-Ming Cheng \> Haizhong Li\\
Department of Mathematics \> Department of Mathematical Sciences\\
Faculty of Science and Engineering \>  Tsinghua University\\
Saga University \> 100084, Beijing\\
840-8502, Saga \> China\\
Japan \> E-mail:hli@math.tsinghua.edu.cn\\
E-mail: cheng@ms.saga-u.ac.jp\\
\end{tabbing}
\end{flushleft}
\begin{flushleft}
\medskip\noindent
\begin{tabbing}
XXXXXXXXXXXXXXXXXXXXXXXXXX*\=\kill Guoxin Wei\\
Department of Mathematics\\
Faculty of Science and Engineering\\
Saga University\\
840-8502, Saga\\
Japan\\
E-mail:wei@ms.saga-u.ac.jp
\end{tabbing}
\end{flushleft}

\end {document}